\documentclass[12pt, reqno]{amsart}
\usepackage{amssymb, amstext, amscd, amsmath}
\usepackage{color}

\usepackage{xy}
\xyoption{all}

\makeatletter
\def\@cite#1#2{{\m@th\upshape\bfseries
[{#1\if@tempswa{\m@th\upshape\mdseries, #2}\fi}]}} \makeatother
\theoremstyle{plain}
\newtheorem{thm}[subsection]{Theorem}
\newtheorem{cor}[subsection]{Corollary}
\newtheorem{prop}[subsection]{Proposition}
\newtheorem{lem}[subsection]{Lemma}

\theoremstyle{definition}
\newtheorem{rem}[subsection]{Remark}

\newtheorem{eg}[subsection]{Example}

\newcommand{\bC}{{\mathbb{C}}}

\newcommand{\bN}{{\mathbb{N}}}

\newcommand{\bQ}{{\mathbb{Q}}}
\newcommand{\bR}{{\mathbb{R}}}

\newcommand{\bT}{{\mathbb{T}}}

\newcommand{\bZ}{{\mathbb{Z}}}

\newcommand{\bn}{{\mathbf{n}}}

\newcommand{\bt}{{\mathbf{t}}}
\newcommand{\bk}{{\mathbf{k}}}
\newcommand{\bx}{{\mathbf{x}}}

  \newcommand{\A}{{\mathcal{A}}}

  \newcommand{\F}{{\mathcal{F}}}

\renewcommand{\O}{{\mathcal{O}}}

\renewcommand{\S}{{\mathcal{S}}}
  
  \newcommand{\U}{{\mathcal{U}}}

  \newcommand{\Z}{{\mathcal{Z}}}

\newcommand{\fD}{{\mathfrak{D}}}

\newcommand{\fF}{{\mathfrak{F}}}

\newcommand{\fM}{{\mathfrak{M}}}

\renewcommand{\phi}{\varphi}
\newcommand{\upchi}{{\raise.35ex\hbox{\ensuremath{\chi}}}}

\newcommand{\qforal}{\quad\text{for all}\quad}

\newcommand{\Ad}{\operatorname{Ad}}

\newcommand{\Aut}{\operatorname{Aut}}

\newcommand{\End}{\operatorname{End}}
\newcommand{\id}{{\operatorname{id}}}
\newcommand{\Inn}{\operatorname{Inn}}

\newcommand{\ran}{\operatorname{Ran}}
\newcommand{\spn}{\operatorname{span}}

\newcommand{\ca}{\mathrm{C}^*}

\newcommand{\Fn}{\mathbb{F}_n^+}
\newcommand{\Fm}{\mathbb{F}_m^+}

\newcommand{\Fth}{\mathbb{F}_\theta^+}

\newcommand{\mt}{\varnothing}

\newcommand{\ol}{\overline}

\newcommand{\y}{{\rm End}(\O_\theta)}
\newcommand{\aut}{{\rm Aut}(\O_\theta)}
\newcommand{\innaut}{{\rm Inn}(\O_\theta)}
\newcommand{\yp}{\lambda_{\epsilon_1}}
\newcommand{\ypp}{\lambda_{\epsilon_2}}
\newcommand{\yppp}{\lambda_{(1,1)}}

\begin{document}
\title[Endomorphisms and Modular Theory]
{Endomorphisms and Modular Theory of 2-Graph C*-Algebras}
\author[D. Yang]
{Dilian Yang$^*$}
\address{Dilian Yang,
Department of Mathematics and Statistics, University of Windsor, Windsor, ON
N9B 3P4, CANADA}
\email{dyang@uwindsor.ca}

\begin{abstract}
In this paper, we initiate the study of endomorphisms and modular theory of the graph C*-algebras $\O_{\theta}$
of a 2-graph $\Fth$ on a single vertex. We prove that
there is a semigroup isomorphism between unital endomorphisms of $\O_{\theta}$ and its
unitary pairs with a \textit{twisted property}. 
We characterize when endomorphisms preserve
the fixed point algebra $\fF$ of the gauge automorphisms and
its canonical masa $\fD$. Some other properties of endomorphisms are also 
investigated. 

As far as the modular theory of $\O_{\theta}$ is concerned,
we show that the algebraic *-algebra generated by the
generators of $\O_{\theta}$ with the inner product induced from a distinguished state $\omega$
is a modular Hilbert algebra. Consequently, we obtain that 
the von Neumann algebra 
$\pi(\O_{\theta})''$ generated by the GNS representation of $\omega$ is 
an AFD factor of type III$_1$,
provided $\frac{\ln m}{\ln n}\not\in\bQ$. 
Here $m,n$ are the numbers of generators of $\Fth$ of degree 
$(1,0)$ and $(0,1)$, respectively. 

This work is a continuation of \cite{DPY1, DPY2} by Davidson-Power-Yang and 
\cite{DY} by Davidson-Yang.
\end{abstract}

\thanks{${}^*$Research partially supported by an NSERC Discovery grant.}

\date{}
\maketitle

\section{Introduction}\label{S:intro}

In 2000, Kumjian-Pask generalized higher rank Cuntz-Kreiger
algebras of Robertson-Steger \cite{RobSte} and introduced
the notion of higher
rank graphs (or $k$-graphs) in \cite{KumPask}. Since then, higher rank graphs
have been attracting a great deal of attention and extensively
studied. See, for example, \cite{FarMuhYee, KP1, KumPask, PaskRRS, P1,
Raeburn, RaeSimYee1, RaeSimYee2, RobSims} and the references therein.

Recently, in \cite{DPY1, DPY2, DY},
Davidson, Power and I have systematically studied an interesting and special class of higher rank
graphs -- 2-graphs with a single vertex, which was initially studied by Power \cite{P1}.
Roughly speaking, those graphs are given concretely in terms of a finite set of generators and
relations of a special type. More precisely, given a permutation $\theta$ of $m\times n$, form a unital
semigroup $\Fth$ with generators $e_1,..., e_m$ and $f_1,...,f_n$ which is free in the
$e_i$'s and free in the $f_j$'s, and has the commutation relations
$e_if_j=f_{j'}e_{i'}$, where $\theta(i, j ) = (i' , j')$
for $1\le i \le m$ and $1\le j \le n$. $\Fth$ is a cancellative semigroup with unique factorization \cite{P1}.

It turns out that 2-graph algebras on a single vertex have very nice structures, and
provide many very interesting and nontrivial phenomena.
We gave a detailed analysis of their representation theory and
completely classified their atomic representations in \cite{DPY1}.
The dilation theory was studied in \cite{DPY2}. Particularly,
it was shown in there that every defect free row
contractive representation of $\Fth$ has a unique minimal *-dilation.
The characterization of the aperiodicity of $\Fth$  and the structure of the graph C*-algebra 
$\O_{\theta}$ were given in \cite{DY}.
From those results, one has a very nice and clear picture of those algebras.
In \cite{DY2}, some of the results in \cite{DPY1, DPY2, DY} were further generalized to $k$-graphs with a
single vertex by Davidson and the author.

The main purpose of this paper is to initiate the study of endomorphisms and modular theory
of 2-graph algebras on a single vertex.
This was motivated by \cite{Cun2} and \cite{CPR}.
This paper can be regarded as a continuation of \cite{DPY1, DPY2, DY}.
In the seminal paper \cite{Cun2}, Cuntz studied the theory of automorphisms of Cuntz algebras $\O_n$.
It is now well-known that there is a
one-to-one correspondence between $\U(\O_n)$,
the  set of unitaries of $\O_n$,
and $\End(\O_n)$, the set of unital endomorphisms of $\O_n$.
We should mention that, very recently, 
the localized automorphisms of Cuntz algebras have been studied in \cite{ConSzy, ConKimSzy, Szy}.

Since the graph C*-algebra $\O_{\theta}$ contains two copies of Cuntz algebras $\O_m$ and $\O_n$, which are
``connected" by the commutation relations of $\Fth$,  one naturally wonders what its unital endomorphisms look like.
In this paper, we show that there is
a semigroup isomorphism between $\End(\O_{\theta})$, the set of unital endomorphisms of $\O_{\theta}$,
and $\U(\O_{\theta})_{\rm twi}^2$, the set consisting of pairs in 
$\U(\O_{\theta})\times \U(\O_{\theta})$ with a \textit{twisted property}.
It is the twisted property that makes an essential difference between the study of
$\End(\O_n)$ and that of $\End(\O_{\theta})$, and makes $\End(\O_{\theta})$ more involved.
However, the twisted property appears here naturally since it decodes the commutation
relations determined by $\theta$.
After obtaining the above isomorphism, we study the
theory of endomorphisms of $\O_{\theta}$ mainly in the vein of \cite{Cun2}.

The second part of this paper is
devoted to investigating the modular theory of $\O_{\theta}$.
This was originally motivated by the index theory of endomorphisms.
There is a rich literature on this topic. See, for example, \cite{CPR, Con} 
and the references therein.
Our first main result here is that
the algebraic *-algebra generated by $\Fth$ with the inner product,
induced from a distinguished faithful state $\omega$ of $\O_{\theta}$, is a modular Hilbert algebra.
We achieve this by obtaining  very explicit expressions of the modular objects
associated with $\omega$ in the celebrated Tomita-Takeski modular theory.
Then we give some partial results on the classification of the von Neumann algebra $\pi(\O_{\theta})''$
generated from the GNS representation of $\omega$.

This paper is organized as follows. The next section provides some preliminaries on 2-graph algebras. In Section \ref{S:End},
we study the general theory of the endomorphisms of $\O_{\theta}$. In particular, we prove that
there is a semigroup isomorphism between unital endomorphisms of $\O_{\theta}$ and
the unitary pairs of $\O_{\theta}$ with a twisted property. Some examples are also given there.
In Section \ref{S:speEnd}, we characterize when endomorphisms or
automorphisms preserve 
the fixed point algebra $\fF$ of the gauge automorphisms and 
its canonical masa $\fD$.
Some other properties of endomorphisms are also 
investigated. 
In Section \ref{S:modular}, the modular theory of 2-graph algebras is
given in detail. We prove that the algebraic *-algebra generated by the
generators of $\O_{\theta}$ with the inner product induced from a distinguished state $\omega$
is a modular Hilbert algebra.
Consequently, we obtain that $\omega$ is a KMS-state with respect to the associated modular automorphism 
group of the von Neumann algebra $\pi(\O_{\theta})''$ generated by $\omega$.
As a further consequence, we show that if $\frac{\ln m}{\ln n}\not\in\bQ$, then
$\pi(\O_{\theta})''$
is an AFD factor of type III$_1$. Those results are given in Section \ref{S:class}.

\section{Preliminaries} \label{S:Pre}

The main source of this section is from \cite{DPY1, DPY2, DY, DY2, KumPask}.

\subsection{2-graphs on a single vertex}
A 2-graph on a single vertex is a unital semigroup
$\Fth$, which is generated by $e_1 , . . . , e_m$  and $f_1 , . . . , f_n$.
The identity is denoted as $\mt$. There are no relations among the $e_i$'s, so they
generate a copy of the free semigroup on $m$ letters, $\Fm$; and there are no
relations on the $f_j$'s, so they generate a copy of $\Fn$. There are commutation
relations between the $e_i$'s and $f_j$'s given by a permutation $\theta$ in $S_{m\times n}$ of $m \times  n$:
\[
e_if_j=f_{j'}e_{i'} \quad \mbox{where}\quad \theta(i, j ) = (i' , j').
\]
The semigroup $\Fth$ has some nice properties. See, for example,  \cite{KP1, KumPask, P1}.
Any word $w\in\Fth$ has fixed numbers of $e$'s and $f$'s
regardless of the factorization. The \textit{degree} of  $w$ is defined as $d(w)=(k,l)$ if there are $k$ $e$'s and
$l$ $f$'s, and the \textit{length} of $w$ is $|w|=k+l$. Moreover, because of the commutation relations, one can
write $w$ according to any prescribed pattern of $e$'s and $f$'s as long as the degree is $(k,l)$.
For instance, we can write $w$ with all $e$'s first, all $f$'s first, or $e$'s and $f$'s alternatively
if $d(w)=(k,k)$.

Recall from \cite{KumPask} that the \textit{graph C*-algebra $\O_{\theta}$} of $\Fth$ is 
the universal $C^*$-algebra generated by
a family of isometries $\{s_u:u\in\Fth\}$ satisfying
$s_{\mt}=I$, $s_{uv}=s_us_v$ for all $u,\ v\in \Fth$,
and the defect free property:
\[
\sum_{i=1}^ms_{e_i}s_{e_i}^*=I=\sum_{j=1}^ns_{f_j}s_{f_j}^*.
\]

It is well-known that $\O_\theta$ has standard generators $s_us_v^*$ 
and $\O_\theta=\overline{\spn}\{s_us_v^*:u,v\in\Fth\}$ (\cite[Lemma 3.1]{KumPask}). 
We extend the degree map $d$ of $\Fth$ to
the generators of $\O_\theta$ as follows: 
$$d(s_us_v^*)=d(u)-d(v)\qforal u,v\in\Fth.$$

\smallskip

To simplify our writing, throughout the paper, 
we use the following multi-index notation:
For all $\bx=(x_1,x_2)\in\bC^2$ with $x_1, x_2\ne 0$  and $\bk=(k_1,k_2)\in\bZ^2$, let
$\bx^\bk:=x_1^{k_1}x_2^{k_2}$. 

\subsection{Gauge automorphisms}

It is well-known that the universal property of $\O_{\theta}$ yields a family of
\textit{gauge automorphisms}
$\gamma_{\bt}$ for $\bt=(t_1, t_2) \in \bT^2$
given by
\[
 \gamma_{\bt}(s_w) =\bt^{d(w)} s_w\qforal w\in\Fth.
\] 
Integrating around $\bT^2$ yields a faithful expectation
\[ \Phi(X) = \int_{\bT^2}  \gamma_{\bt}(X) \,d\bt\]
onto the fixed point algebra $\O_{\theta}^\gamma$ of $\gamma$.

It turns out that
$$\fF:=\O_{\theta}^\gamma=\ol{\bigcup_{k\ge 1} \fF_k}$$
is an $(mn)^\infty$-UHF algebra,
where
$\fF_k=\ol\spn\{ s_us_v^* : d(u)=d(v)=(k,k)\}$ ($k\in\bN$)
is the full matrix algebra $M_{(mn)^k}$.

For each $\bn=(n_1,n_2)\in\bZ^2$, define a mapping $\Phi_\bn$ on $\O_{\theta}$ via
\[
\Phi_\bn(X)=\int_{\bT^2}\bt^{-\bn}\gamma_\bt(X) d\bt \quad \mbox{for all}\quad X\in\O_{\theta}.
\]
Then $\Phi_\bn$ acts on generators via
$$
\Phi_\bn(X)=\left\{
\begin{array}{ll}
X, & \quad \mbox{if}\quad d(X)=\bn,\\
0, &\quad \mbox{otherwise}.
\end{array}\right.
$$
So the fixed point algebra $\fF$ is nothing but
$\ran\Phi_\mathbf{0}$, which is spanned by the words of
degree $(0,0)$. 
For every $X\in\O_{\theta}$, we also have the formal series
$
X\sim \sum_{\bn\in\bZ^2}\Phi_\bn(X)
$
with a Cesaro convergence of the series.  Refer to \cite{ HJP} for the details.

Let
$$\fD=\ol\spn\{s_ws_w^*: d(w)=(k,k),\ k\in\bN\},$$
the canonical masa in $\fF$,
and
$$N(\fD)=\{U\in\U(\O_{\theta}): U\fD U^*\subseteq \fD\},$$
the (unitary) normalizer of $\fD$.

Notice that, in the literature, it is usual to interpret
$\fF$ (resp. $\fD$) as the C*-algebras generated by
$\{s_us_v^* : d(u)=d(v)=(k,l), k, l\in\bN\}$
(resp.
$\{s_ws_w^*: d(w)=(k,l),\ k,l\in\bN\}$).
But they are the same as those given above
because of the defect free property
$\sum_{i=1}^ms_{e_i}s_{e_i}^*=I=\sum_{j=1}^ns_{f_j}s_{f_j}^*$.
Refer to \cite{DY} for more details.
We find that it is more convenient to use the above way for our purpose.
For example, when we consider the restriction $\lambda_{(U,V)}|_\fF$
of an endomorphism $\lambda_{(U,V)}$ determined by $(U,V)$,
it suffices to take care of the generators $X=s_us_v^*$
with $d(u)=d(v)=(k,k)$.
Using the commutation relations, one can rewrite $X$ as 
$X=s_{u_1\cdots u_k}s_{v_1\cdots v_k}^*$,
where $u_i,v_i\in\Fth$ with $d(u_i)=d(v_i)=(1,1)$ ($i=1,...,k$).
Then the action of $\lambda_{(U,V)}$ on $X$ is now completely determined by the unitary
$W=U\yp(V)=V\ypp(U)$ (see Section \ref{S:End} for the notation).
We will often use this simple and useful observation later.

\subsection{Notation and conventions}
We end this section with introducing some notation and conventions.
Let $\bZ_+$ denote the set of all non-negative integers.
If $\A$ is a unital C*-algebra, by $\U(\A)$ we mean the set of all unitaries of $\A$.
Let $\End(\A)$ denote the semigroup of unital endomorphisms of $\A$, and
$\rm Aut(\A)$, $\Inn(\A)$ the groups of automorphisms, inner automorphisms of $\A$, respectively.
The notation  $\Z(\A)$ stands for the center of $\A$. In this paper, by endomorphisms we always mean
\textit{unital} endomorphisms.

\section{Endomorphisms of 2-graph algebras}\label{S:End}

In this section, we will study the general theory of unital endomorphisms of $\O_\theta$.
Some examples will also be  given.

\subsection{General Theory}\label{SubS: general}

Let $\Fth$ be a 2-graph on a single vertex generated by $e_1,\ldots, e_m$ and 
$f_1,\ldots, f_n$. It is easy to see that, for $(p,q)\in\bZ_+^2$,
the completely positive maps $\lambda_{(p,q)}$ on $\O_\theta$ given by
\begin{eqnarray*}
\lambda_{(p,q)}(X)=\sum_{d(w)=(p,q)} s_wXs_w^*\quad \mbox{for all}\quad X\in\O_\theta
\end{eqnarray*}
are endomorphisms of $\O_\theta$. Those endomorphisms are said to be \textit{canonical}.
Let $\epsilon_1$ and $\epsilon_2$ be the standard generators of $\bZ_+^2$. One can easily check that
$\lambda_{(p,q)}=\yp^p\ypp^q=\ypp^q\yp^p$, and
that $\lambda_{(p,q)}(X)s_w=s_wX$ for all $X\in \O_\theta$ and $w$ with $d(w)=(p,q)$.

Let $\U(\O_\theta)_{\rm twi}^2\subseteq\U(\O_\theta)\times \U(\O_\theta)$ be the family of unitary pairs
satisfying a \textit{twisted property}. More precisely,
$$\U(\O_\theta)_{\rm twi}^2=\{(U,V)\in\U(\O_\theta)\times \U(\O_\theta): \ U\yp(V)=V\ypp(U)\}.$$
Then there is a bijective correspondence between $\U(\O_\theta)_{\rm twi}^2$ and $\y$, as shown below.

\begin{thm}\label{T:End}
The mapping
$\Psi: \U(\O_\theta)_{\rm twi}^2\to \y,$
$$(U,V)\mapsto\lambda_{(U,V)}:\ s_{e_i}\mapsto Us_{e_i},\ s_{f_j}\mapsto Vs_{f_j}\ (i=1,...,m,\  j=1,...,n)$$
is a bijection.
Its inverse is given by 
$$\lambda\mapsto\left(\sum_{i=1}^m\lambda(s_{e_i}) s_{e_i}^*,\ \sum_{j=1}^n\lambda(s_{f_j}) s_{f_j}^*\right).$$
\end{thm}

\begin{proof}
Let $(U,V)\in \U(\O_\theta)_{\rm twi}^2$. Define
$$\lambda_{(U,V)}(s_{e_i})=Us_{e_i}, \ \lambda_{(U,V)}(s_{f_j})=Vs_{f_j},\ i=1,...,m,\ j=1,...,n.$$
In what follows, we prove that 
$\lambda_{(U,V)}$ preserves the commutation relations given by $\theta$. That is,
if $e_if_j=f_{j'}e_{i'}$ where $\theta(i,j)=(i',j')$, then
$\lambda_{(U,V)}(s_{e_i})\lambda_{(U,V)}(s_{f_j})
 =\lambda_{(U,V)}(s_{f_{j'}})\lambda_{(U,V)}(s_{e_{i'}})$.
Indeed, let
$W:=U\yp(V)=V\ypp(U)$.
Since $e_if_j=f_{j'}e_{i'}$, we have
\begin{alignat*}{2}
U\sum_{i=1}^ms_{e_i}Vs_{e_i}^*=W & \Longrightarrow Us_{e_i}V=Ws_{e_i} \Longrightarrow Us_{e_i}Vs_{f_j}=Ws_{e_if_j};\\
V\sum_{j'=1}^ns_{f_{j'}}Us_{f_{j'}}^*=W & 
\Longrightarrow Vs_{f_{j'}}U=Ws_{f_{j'}} \Longrightarrow Vs_{f_{j'}}Us_{e_{i'}}=Ws_{f_{j'}e_{i'}}.
\end{alignat*}
Thus
$Us_{e_i}Vs_{f_j}=Vs_{f_{j'}}Us_{e_{i'}}$,
which says that $\lambda_{(U,V)}$ preserves the commutation relations.

By the universal property of $\O_\theta$, the mapping $\lambda_{(U,V)}$ 
defines an endomorphism of $\O_\theta$.

Conversely, suppose $\lambda\in\End(\O_\theta)$. Set
$$
U=\sum_{i=1}^m\lambda(s_{e_i}) s_{e_i}^*,\quad
V=\sum_{j=1}^n\lambda(s_{f_j}) s_{f_j}^*,\quad
W=\sum_{i=1}^m\sum_{j=1}^n\lambda(s_{e_if_j}) s_{e_if_j}^*.
$$
Then it is straightforward to verify that $U, V, W$ all are unitaries, and that
\begin{align*}
\lambda(s_{e_i})=Us_{e_i},\quad
\lambda(s_{f_j})=Vs_{f_j},\quad
\lambda(s_{e_if_j})=Ws_{e_if_j}.
\end{align*}
Furthermore, $U, V, W$ have
the following relations:
\begin{align*}
U^*W&=\sum_{k=1}^ms_{e_k}\lambda(s_{e_k})^* \sum_{k,l}\lambda(s_{e_kf_l}) s_{e_kf_l}^*\\
    &=\sum_{k=1}^ms_{e_k} \left(\sum_{l=1}^n\lambda(s_{f_l}) s_{f_l}^*\right)s_{e_k}^*\\
    &=\sum_{k=1}^ms_{e_k} V s_{e_k}^*\\
    &=\yp(V).
\end{align*}
Similarly,
$$V^*W=\ypp(U).$$
Therefore, $(U,V)$ satisfies the twisted property required in the proposition:
$U\yp(V)=V\ypp(U)(=W)$. We proved $(U,V)\in \U(\O_\theta)_{\rm twi}^2$.

It is easy to check that the above two processes are inverses of each other. This ends the proof. 
\end{proof}

Some remarks are in order.

\begin{rem}\label{R:canend}
Because of Theorem \ref{T:End}, we can and do use
$\lambda_{(U,V)}$ to denote the endomorphism of $\O_\theta$
uniquely determined by the unitary pair $(U,V)\in\U(\O_\theta)_{\rm twi}^2$.
Using this notation, the canonical endomorphism $\lambda_{(p,q)}$ is
equal to $\lambda_{(U,V)}$, where
\[U=\sum_{i=1}^m\sum_{d(w)=(p,q)}s_{we_i}s_{e_iw}^*, \quad
  V=\sum_{j=1}^n\sum_{d(w)=(p,q)}s_{wf_j}s_{f_jw}^*.\]

\end{rem}

\begin{rem}
From  the proof of Theorem \ref{T:End}, we actually have
$W=U\yp(V)=V\ypp(U)$.
So any two of the three unitaries $U,V,W$
completely determine the endomorphism $\lambda_{(U,V)}$.
\end{rem}

\begin{rem}
As we shall see later, the twisted property in Theorem \ref{T:End} makes the study of 
$\End(\O_\theta)$ more interesting, and gives the essential difference between the studies of
$\End(\O_\theta)$ and $\End(\O_n)$. 
\end{rem}

We should also admit that, in general,  it is not easy
to check if the twisted property in Theorem \ref{T:End} holds for a given unitary pair.
However, we do have a lot of examples.  Before giving our examples, let us 
first consider the composition of two endomorphisms.  
We will see that  $\U(\O_\theta)_{\rm twi}^2$ is a untial semigroup
with the multiplication induced from the composition. 
The following lemma is adapted from \cite[Proposition 1.1]{Cun2}
and the discussion immediately following. 

\begin{lem}\label{L:mult}
Endomorphisms of $\O_\theta$ have the following properties. 
\begin{itemize}
\item[(i)] If $\lambda_{(U_i,V_i)}\in\End(\O_\theta)$ ($i=1,2$), then
$\lambda_{(U_2,V_2)}\lambda_{(U_1,V_1)}=\lambda_{(U,V)}$,
where
$$U=\lambda_{(U_2,V_2)}(U_1)U_2,\quad V=\lambda_{(U_2,V_2)}(V_1)V_2.$$
\item[(ii)] Suppose $\Fth$ is aperiodic and $\lambda_{(U,V)}\in\End(\O_\theta)$.
Then
$$
\lambda_{(U,V)}\in\aut\Longleftrightarrow
\lambda_{(U,V)}(U_0)=U^*,\ \lambda_{(U,V)}(V_0)=V^*
$$
for some
$U_0,V_0\in\U(\O_\theta)$.
\item[(iii)] The mapping
\begin{alignat*}{2}
\phi : \U(\O_\theta) \to \innaut, \quad
 W \mapsto \lambda_{(U,V)}
\end{alignat*}
is surjective, where $(U,V)=(W\yp(W)^*, W\ypp(W)^*)$. Moreover,
$\ker\phi\subseteq\Z(\O_\theta)$.
\end{itemize}
\end{lem}

\begin{proof}
(i)
For simplicity, put $\lambda_i:=\lambda_{(U_i,V_i)}$ ($i=1,2$). Then
\begin{align*}
\lambda_2(\lambda_1(s_{e_i}))&=\lambda_2(U_1s_{e_i})
=\lambda_2(U_1)\lambda_2(s_{e_i})=\lambda_2(U_1)U_2s_{e_i}=Us_{e_i},\\
\lambda_2(\lambda_1(s_{f_j}))&=\lambda_2(V_1s_{f_j})
=\lambda_2(V_1)\lambda_2(s_{f_j})=\lambda_2(V_1)V_2s_{f_j}=Vs_{f_j},
\end{align*}
where $U=\lambda_2(U_1)U_2, \ V=\lambda_2(V_1)V_2$.
Hence the composition endomorphism $\lambda_2\lambda_1$ is
determined by $(U,V)$, i.e., $\lambda_2\lambda_1=\lambda_{(U,V)}$.

(ii) It suffices to check the sufficiency. Since $\Fth$ is
aperiodic, $\O_\theta$ is simple. See, e.g., \cite{DY} (or Lemma \ref{L:masa}).
So
$\lambda_{(U,V)}$ is injective. But if
$\lambda_{(U,V)}(U_0)=U^*$ and $\lambda_{(U,V)}(V_o)=V^*$ for some
$U_0,V_0\in\U(\O_\theta)$, then
$\lambda_{(U,V)}(U_0s_{e_i})=\lambda_{(U,V)}(U_0)\lambda_{(U,V)}(s_{e_i})=U^*Us_{e_i}=s_{e_i}$, 
and,
similarly, $\lambda_{(U,V)}(V_0s_{f_j})=s_{f_j}$. This gives the surjectivity of
$\lambda_{(U,V)}$. Hence $\lambda_{(U,V)}$ is an automorphism.

(iii) It is sufficient to show that $\lambda_{(U,V)}\in\Inn(\O_\theta)$ iff $U=W\yp(W^*)$ and $V=W\ypp(W^*)$
for some $W\in\U(\O_\theta)$. Clearly,
$\lambda_{(U,V)}\in\Inn(\O_\theta)$  
iff $\lambda_{(U,V)}(s_{e_i})=Ws_{e_i}W^*$ 
and 
$\lambda_{(U,V)}(s_{f_j})=Ws_{f_j}W^*$ ($i=1,...,m,j=1,...,n$)
for some $W\in\U(\O_\theta)$.
But we also have $\lambda_{(U,V)}(s_{e_i})=Us_{e_i}$ and
$\lambda_{(U,V)}(s_{f_j})=Vs_{f_j}$.
Thus
\begin{align*}
&\lambda_{(U,V)}\in\Inn(\O_\theta)\\
&\Leftrightarrow
Us_{e_i}=Ws_{e_i}W^*, \quad  Vs_{f_j}=Ws_{f_j}W^*\ 
                  \text{for some}\  W\in\U(\O_\theta)\\
&\Leftrightarrow
Us_{e_i}s_{e_i}^*=Ws_{e_i}W^*s_{e_i}^*,\quad  
Vs_{f_j}s_{f_j}^*=Ws_{f_j}W^*s_{f_j}^*\\
&\Leftrightarrow
\sum_{i=1}^mUs_{e_i}s_{e_i}^*=\sum_{i=1}^mWs_{e_i}W^*s_{e_i}^*,\quad 
\sum_{j=1}^nVs_{f_j}s_{f_j}^*=\sum_{j=1}^nWs_{f_j}W^*s_{f_j}^*\\
&\Leftrightarrow
U=W\yp(W^*),\quad V=W\ypp(W^*),
\end{align*}
where $i=1,...,m, \ j=1,...,n$. 
(The property of
$(U,V)\in\U(\O_\theta)_{\rm twi}^2$ is automatic.) This ends the proof.
\end{proof}

By Theorem \ref{T:End} and Lemma \ref{L:mult} (i), we immediately have

\begin{cor}
Let $(U_1,V_1), (U_2,V_2))\in\U(\O_\theta)_{\rm twi}^2$. Then,
under the multiplication
$$(U_2,V_2)\cdot (U_1,V_1)
:=(\lambda_{(U_2,V_2)}(U_1)U_2,\lambda_{(U_2,V_2)}(V_1)V_2),$$
the set  $\U(\O_\theta)_{\rm twi}^2$ is a unital semigroup, and the bijection
$\Psi$ in Theorem \ref{T:End} is a unital semigroup isomorphism.
\end{cor}

\subsection{Examples}\label{SubS:eg}

We now give some examples of non-canonical endomorphisms of $\O_{\theta}$. 

\begin{eg}
Consider $\Fth$ with $m=n$ and $\theta$ the flip relation: $e_if_j=f_ie_j$ for all $i,j=1,...,m$.
Let $U$ be any unitary of $\O_\theta$. Then $(U,U)\in\U(\O_\theta)_{\rm twi}^2$ and so it gives an
endomorphism of $\O_\theta$ by Theorem \ref{T:End}.
Indeed,
\begin{align*}
(U,U)\in\U(\O_\theta)_{\rm twi}^2
&\Leftrightarrow
U\yp(U)=U\ypp(U)\Leftrightarrow \yp(U)=\ypp(U)\\
&\Leftrightarrow
(\sum_{k=1}^ms_{e_k}Us_{e_k}^*)s_{e_i}=(\sum_{l=1}^ms_{f_l}Us_{f_l}^*)s_{e_i}
\ (i=1,...,m)\\
&\Leftrightarrow
s_{f_j}^*s_{e_i}U=s_{f_j}^*(\sum_{l=1}^ms_{f_l}Us_{f_l}^*)s_{e_i} \ (i,j=1,...,m)\\
&\Leftrightarrow
s_{f_j}^*s_{e_i}U=Us_{f_j}^*s_{e_i}\ (i,j=1,...,m).
\end{align*}
Therefore, $(U,U)\in\U(\O_\theta)_{\rm twi}^2$ iff $U$ commutes with $s_{f_j}^*s_{e_i}$ $(i,j=1,...,m$).
But for flip algebras, we have $s_{f_j}^*s_{e_i}=\delta_{i,j}s_{f_i}^*s_{e_i}$, and
$s_{f_i}^*s_{e_i}$ is a unitary in $\Z(\O_\theta)$ (cf. \cite{DY}).
Here, as usual, $\delta_{i,j}$ denotes the Kronecker delta function.
Hence we have $(U,U)\in\U(\O_\theta)_{\rm twi}^2$ for any $U\in\U(\O_\theta)$.

Furthermore, one can check that if $U=W\lambda_{\epsilon_{1}}(W)^{*}$ for some $W\in\U(\O_{\theta})$ then
$\lambda_{(U,U)}$ is an inner automorphism, as
$ W\lambda_{\epsilon_{1}}(W)^{*}=W\lambda_{\epsilon_{2}}(W)^{*}$
in this case. 

\end{eg}

\begin{eg}\label{Eg:center}
Let $\Fth$ be a 2-graph.
Let $U,V\in\U(\O_\theta)$. If 
$$UV=VU, \ Us_{f_j}=s_{f_j}U, \ Vs_{e_i}=s_{e_i}V\  (i=1,...,m,\ j=1,...,n),$$ 
then one can easily check that $(U,V)\in\U(\O_\theta)_{\rm twi}^2$, and so $(U,V)$
determines an endomorphism. In particular, $(U,V)$ with $U,V
\in\Z(\O_\theta)$ gives an endomorphism of $\O_\theta$.
\end{eg}

\begin{eg}\label{Eg:id}
Consider $\Fth$ with the identity relation $\theta$: $e_if_j=f_je_i$ for all $i=1,...,m,\  j=1,...,n$.
Then we have $\O_\theta\cong\O_m\otimes \O_n$ by \cite[Corollary 3.5 (iv)]{KumPask}.
Let $U$ (resp. $V$) be unitaries in the Cuntz algebras 
$\ca(s_{e_1},...,s_{e_m})\cong\O_m$ 
(resp. $\ca(s_{f_1},...,s_{f_n})\cong\O_n$). Then, from Example \ref{Eg:center},
$(U,V)$ gives an endomorphism of $\O_\theta$. Furthermore, if $U$ (resp. $V$)
determines an automorphism of $\O_m$ (resp. $\O_n$), then $\lambda_{(U,V)}\in\aut$
by Lemma \ref{L:mult} (ii). Moreover, it is easy to see that 
$\lambda_{(U,V)}=\lambda_U\otimes \lambda_V$
by considering the actions on generators,  where 
$\lambda_U$ (resp. $\lambda_V$) are endomorphisms of $\O_m$ (resp. $\O_n$)  
as defined in \cite{Cun2}. 
\end{eg}

In \cite{Arc}, it was proved that the ``flip-flop'' unitary
$U=s_{e_2}s_{e_1}^*+s_{e_1}s_{e_2}^*$
gives an outer automorphism of $\ca(s_{e_1},s_{e_2})\cong \O_2$.
In the following example, we give an outer automorphism of $\O_\theta$
 by using two ``mixing'' flip-flop unitaries.

\begin{eg}\label{Eg:flipalg}
Suppose $m=n$ and $\theta$ is the flip relation. Let
$$
U=\sum_{j=1}^m s_{s_{f_j}}s_{e_j}^*\ (=s_{e_i}^*s_{f_i}).
$$
Notice that $U\in\Z(\O_\theta)$ (cf. \cite{DY}). From Example
\ref{Eg:center}, $(U,U^*)$ gives an endomorphism $\lambda_{(U,U^*)}$.
A simple calculation shows that $\lambda_{(U,U^*)}$ is an involution:
$\lambda_{(U,U^*)}^2=\id$. In particular, $\lambda_{(U,U^*)}$ is an automorphism.

In what follows, we show that $\lambda_{(U,U^*)}$ is actually
outer. To the contrary, suppose that
$\lambda_{(U,U^*)}=\Ad W$ for some $W\in\U(\O_\theta)$. Then
\begin{alignat*}{2}
\lambda(s_{e_i})&=Us_{e_i}=Ws_{e_i}W^*=s_{f_i},\\
\lambda(s_{f_i})&=U^*s_{f_i}=Ws_{f_i}W^*=s_{e_i}
\end{alignat*}
for all $i=1,...,m$. Therefore
$$Ws_{e_i}W^*=s_{f_i}=W^*s_{e_i}W,\quad Ws_{f_i}W^*=s_{e_i}=W^*s_{f_i}W.$$
We have
$W^2s_{e_i}=s_{e_i}W^2$ and $W^2s_{f_i}=s_{f_i}W^2$ for $i=1,...,m$. So
$W^2\in\Z(\O_\theta)$. Also, from above we have $U=W\yp(W)^*$. This implies
$U^*W=\yp(W)$ and  $WU=W^2\yp(W)^*=\yp(W)$ as $W^2\in\Z(\O_\theta)$.
So $U^*W=WU$. Since $U\in\Z(\O_\theta)$, we have
$UW=U^*W$. Namely, $U=U^*$. This is ridiculous.
(Indeed, $U=U^*\Rightarrow s_{e_i}^*s_{f_i}s_{f_j}=s_{f_i}^*s_{e_i}s_{f_j}
\Rightarrow s_{e_i}^*s_{f_if_j}=s_{e_j}
\Rightarrow s_{e_ie_j}^*s_{f_if_j}=I \Rightarrow s_{e_ie_j}=s_{f_if_j}$, a contradiction.)
Therefore
$\lambda_{(U,V)}$ is outer.

\smallskip

By \cite{DPY2, DY}, we have $\O_\theta\cong C(\bT)\otimes \O_m$.
From this point of view, one can see that the above
$\lambda_{(U,V)}\cong \phi\otimes \id$, where $\phi\in\rm
Aut(C(\bT))$ given by $\phi(z)=\bar z$ (the conjugate of $z$),
and $\id$ is the trivial
automorphism of $\O_m$.

\medskip
Curiously, if we replace $(U,V)$ in Example \ref{Eg:flipalg} to 
$(U,V)=(s_{e_2}s_{e_1}^*+s_{e_1}s_{e_2}^*, s_{f_2}s_{f_1}^*+s_{f_1}s_{f_2}^*)$, it is rather easy to 
see that $\lambda_{(U,V)}$ is also an automorphism. But we do not know if it is outer in this case.  
The essential difference is that $U$ is not in $\Z(\O_\theta)$ any more.

It follows from Example \ref{Eg:center} that each pair
$(U^k,{U^*}^k)$ with $k\ge 1$ indeed gives an endomorphism $\lambda_k:=\lambda_{(U^k,{U^*}^k)}$.
It is not hard to see that $\lambda_k$ with $k\ge 2$ is not an
automorphism, and that $\lambda_k|_\fF=\id$.
Since we will not use this fact later, we omit the details here.
\end{eg}

\begin{eg}
Let $m=n$ and $\theta$ be the identity relation. Let
$
U=\sum_{j=1}^m s_{f_j}s_{e_j}^*.
$
One can easily check that $(U,U^*)\in\U(\O_\theta)_{\rm twi}^2$, and so it determines an automorphism $\lambda_{(U,U^*)}$.
But, unlike Example \ref{Eg:flipalg}, we do not know if $\lambda_{(U,U^*)}$ is  outer in this case.
Using the same argument as in Example \ref{Eg:flipalg}, we still have
$W^2\in\Z(\O_\theta)$, and so $W^2$ is a scalar as $\O_\theta$ is now simple \cite{DY}.
However, the issue here is that $U$ is not in $\Z(\O_\theta)$ any more.
\end{eg}

\section{Endomorphisms preserving some Subalgebras}  \label{S:speEnd}

In this section, we study endomorphisms of $\O_\theta$ which preserve 
the fixed point algebra $\fF$ of the gauge automorphisms $\gamma_\bt$ $(\bt\in\bT^2$)
and its diagonal subalgebra $\fD$. 
Recall from Section \ref{S:Pre} that $\fD$ is the canonical masa of $\fF$.
When $\Fth$ is aperiodic, we actually have more:
$\fD$ is also a masa in $\O_\theta$, and $\fF$ has the trivial relative commutant;
moreover, the converse is also true.
We should mention that it is well-known that these hold true for Cuntz algebras $\O_n$
(see, e.g., \cite[Section 1] {Cun2} and \cite[Corollary 3.3]{DR}).
\begin{lem}
\label{L:masa}
The following statements are equivalent.
\begin{itemize}
\item[(i)] $\Fth$ is aperiodic.
\item[(ii)] $\O_\theta$ is simple.
\item[(iii)] $\fD'\cap\O_\theta=\fD$.
\item[(iv)] $\fF'\cap\O_\theta=\bC I$.
\end{itemize}
\end{lem}

\begin{proof}
By \cite{DY} or \cite{RobSims}, (i) and (ii) are equivalent.

That (i) is equivalent to (iii) is the main result of \cite{Hope}. (Notice that
$\fD$ here is the same as  $\fD$ in \cite{Hope} in terms
of groupoid terminology.)

(iii)$\Rightarrow$(iv): This is directly from the simple fact that the relative commutant
of an algebra is contained in that of a subalgebra. 

(iv)$\Rightarrow$(i):
If $\Fth$ is periodic, by a result in \cite{DY}, we have
$\rm C(\bT)\subseteq\Z(\O_\theta)\subseteq\fF'\cap\O_\theta$ and so $\fF'\cap\O_\theta\ne \bC I$.
\end{proof}

Because the proofs in the rest of this subsection need the property 
that $\fD$ is a masa in $\O_{\theta}$, we assume that $\Fth$ is aperiodic 
by Lemma  \ref{L:masa}.

The proofs of Propositions \ref{P:invF} and \ref{P:invD} below
are adapted from those of \cite[Proposition 1.2]{Cun2} and  \cite[Propositions 1.3--1.5]{Cun2},
respectively. 
We sketch the proofs here, and only give the full details which are different from there. 

\begin{prop}\label{P:invF}
Suppose that $\Fth$ is aperiodic.
Let $\lambda_{(U,V)}\in\y$ and $W:=U\yp(V)(=V\ypp(U))$.
Then we have the following.
\begin{itemize}
\item[(i)] $W\in\fF \Rightarrow \lambda_{(U,V)}(\fF)\subseteq \fF$. Moreover, 
$\lambda_{(U,V)}(\fF)\subseteq \fF\Rightarrow W\in\fF$ provided that 
$\lambda_{(U,V)}(\fF)'\cap \O_\theta=\bC I$.
\item[(ii)] If $U,V\in\fF$, then
$\lambda_{(U,V)}\in\aut\Longleftrightarrow \lambda_{(U,V)}(\fF)=\fF$.
\end{itemize}
\end{prop}

\begin{proof}
(i)
We first claim that
\begin{eqnarray}\label{E:limInv}
\lambda_{(U,V)}(X)=\lim_{k\to\infty}\Ad(W\yppp(W)\cdots \yppp^k(W))(X)\quad\text{for all}\quad X\in\fF.
\end{eqnarray}
It suffices to check \eqref{E:limInv} for the generators $X$ of $\fF$.
But if $X_1\in\fF_1$, then $X_1$ can be written as $X_1=s_{u_1}s_{v_1}^*$
with $d(u_1)=d(v_1)=(1,1)$. Recall that $\lambda_{(U,V)}(s_{w})=Ws_{w}$
for all $w\in\Fth$ with $d(w)=(1,1)$. 
So we have
$\lambda_{(U,V)}(X_1)=\Ad (W)(X_1)$.
Now one can prove inductively 
$$\lambda_{(U,V)}(X_{k})
=\Ad(W\yppp(W)\cdots \yppp^{k-1}(W))(X_k) \quad \mbox{for all}\quad X_k\in \fF_k.$$
To this end, first notice that simple calculations yield the following relations:
For all $w\in\Fth$ with $d(w)=(1,1)$, 
$$s_w\lambda_{(1,1)}^{k-1}(W)
=\lambda_{(1,1)}^{ k}(W)s_w\ (i=1,..., m, j=1,..., n, k\ge 1).$$
Let 
$X_{ k}
=s_{u_1\cdots u_{ k}}s_{v_1\cdots v_k}^*\in\fF_k$,
where $d(u_i)=d(v_i)=(1,1)$ ($i=1,..., k$).  
Put $X_{ k-1}=s_{u_2\cdots u_k} s_{v_2\cdots v_k}^*$. Obviously, $X_{k-1}\in\fF_{ k-1}$. 
Then use the inductive assumption and the relations given above to obtain   
\begin{align*}
\lambda_{(U,V)}(X_{k})
&=\lambda_{(U,V)}(s_{u_1}X_{ k -1}s_{v_1}^*)\\
&=Ws_{u_1}\Ad(W\yppp(W)\cdots \yppp^{ k-2}(W))(X_{k-1})s_{v_1}^*W^*\\
&=W\cdot s_{u_1}W\cdot \yppp(W)\cdots \yppp^{ k-2}(W)X_{k-1}\cdot \\
&\quad\hskip .2cm  \yppp^{k-2}(W)^*\cdots \yppp(W)^*\cdot W^* s_{v_1}^*\cdot W^*\\
&=W\cdot \yppp(W)s_{u_1}\cdot \yppp(W)\cdots \yppp^{ k-2}(W)X_{k-1}\cdot\\
&\quad\hskip .2cm \yppp^{ k-2}(W)^*\cdots \yppp(W)^*\cdot s_{v_1}^*\yppp(W)^*\cdot W^*\\
&=\cdots\\
&=W\yppp(W)\cdots \yppp^{k-1}(W) s_{u_1}X_{k-1}s_{v_1}^*
      \yppp^{ k-1}(W)^*\cdots \yppp(W)^*W^* \\
&=\Ad(W\yppp(W)\cdots \yppp^{ k-1}(W))(X_{k}).
\end{align*} 
This ends the proof of our claim.

That $W\in\fF \Rightarrow \lambda_{(U,V)}(\fF)$ is directly derived from \eqref{E:limInv}.

We now assume that $\lambda_{(U,V)}(\fF)'\cap \O_\theta=\bC I$.
It is easy to check
\begin{alignat}{1}
\label{E:lambda}
\gamma_\bt\lambda_{(U,V)}\gamma_{\bt}^{-1}=\lambda_{(\gamma_\bt(U),\gamma_\bt(V))}.
\end{alignat}
Since $\lambda_{(U,V)}(\fF)\subseteq\fF$ and $\gamma_\bt(X)=X$ for all $X\in\fF$,
then \eqref{E:lambda} yields
\begin{eqnarray}
\label{E:eqF}
\lambda_{(U,V)}(X)=\lambda_{(\gamma_\bt(U),\gamma_\bt(V))}(X) \quad \mbox{for all}\quad  X\in\fF.
\end{eqnarray}
Direct calculations give
\begin{eqnarray}
\label{E:W}
\gamma_\bt(U)\yp(\gamma_\bt(V))=\gamma_\bt(U\yp(V))=\gamma_\bt(W).
\end{eqnarray}
We now claim $W^*\gamma_\bt(W)\in\lambda_{(U,V)}(\fF)'\cap\O_\theta$.
Indeed, if $X_1=s_{u_1}s_{v_1}^*$ with $d(u_1)=d(v_1)=(1,1)$, 
from \eqref{E:eqF} and \eqref{E:W}
we derive
$WX_1W^*=\gamma_\bt(W)X_1\gamma_\bt(W)^*$.
Thus
$$W^*\gamma_\bt(W) \tilde X_1=\tilde X_1W^*\gamma_\bt(W),\quad \text{where}\quad \tilde X_1:=X_1.$$
Then putting
$X_2=s_{u_1u_2}s_{v_1v_2}^*$
to \eqref{E:eqF} and using \eqref{E:W}, we have 
$$
W^*\gamma_\bt(W)\tilde X_2=\tilde X_2W^*\gamma_\bt(W),
$$
where 
$\tilde X_2:=s_{u_1}Ws_{u_2}s_{v_2}^*W^*s_{v_1}^*$.
Set 
$$
\tilde X_n:=s_{u_1}Ws_{u_2}\cdots W s_{u_n} s_{v_n}^*W^* \cdots s_{v_2}^*W^*s_{v_1}^*.
$$
Here the degrees of all above $u_i$'s and $v_i$'s are $(1,1)$.
By induction, we can obtain 
$$
W^*\gamma_\bt(W) \tilde X_n=\tilde X_n W^*\gamma_\bt(W).
$$
Furthermore, from the above process one also has
$\tilde X_n\in\fF$ as $\lambda_{(U,V)}(\fF)\subseteq\fF$.
Therefore, we have
$W^*\gamma_\bt(W)$ commutes with every element in
$$
\tilde \fF=\ol{\bigcup_{n\ge 1}\spn\{\tilde X_n\text{'s}\}}= \overline{\spn}\{s_u\lambda_{(U,V)}(\fF)s_v^*: d(u)=d(v)=(1,1)\}.
$$
Now one can check that 
$\gamma_\bt(W)W^*\in \lambda_{(U,V)}(\fF)'\cap\O_\theta$.
This ends the proof of our claim.

Therefore, from our assumption, we obtain
$\lambda_\bt(W)=\alpha W$ for some $\alpha\in\bT$.
This forces $W\in\fF$. We are done.

(ii) Let $U,V\in\fF$. Suppose $\lambda_{(U,V)}\in\aut$.
Obviously, $W\in\fF$. This implies $\lambda_{(U,V)}(\fF)\subseteq\fF$ from
\eqref{E:limInv}.
Also, as $\gamma_\bt(U)=U$ and $\gamma_\bt(V)=V$,
the identity \eqref{E:lambda} yields
$\lambda_{(U,V)}^{-1}\gamma_{\bt}=\gamma_\bt\lambda_{(U,V)}^{-1}$.
Thus, $\lambda_{(U,V)}^{-1}(\fF)\subseteq\fF$.
Therefore, $\lambda_{(U,V)}(\fF)=\fF$.

For the other direction, since $\Fth$ is aperiodic, $\O_\theta$
is simple by Lemma \ref{L:masa}. So $\lambda_{(U,V)}$ is injective. Since $U^*,V^*\in\fF$ and
$\lambda_{(U,V)}(\fF)=\fF$, there are
$U_0,V_0\in\fF$ such
that $\lambda_{(U,V)}(U_0)=U^*$, $\lambda_{(U,V)}(V_0)=V^*$. This implies that $\lambda_{(U,V)}$ is surjective
as in the proof of Lemma \ref{L:mult} (ii).
So $\lambda_{(U,V)}$ is an automorphism.
\end{proof}

Clearly, by Remark \ref{R:canend}, all canonical endomorphisms $\lambda_{(p,q)}$ satisfy the corresponding unitary pairs
$(U,V)\in\fF\times \fF$. So $\lambda_{(p,q)}(\fF)\subseteq \fF$, and
$\lambda_{(p,q)}\in\Aut(\O_\theta)$ if and only if $\lambda_{(p,q)}(\fF)=\fF$.

Before stating the following result, we recall that $N(\fD)$ is the unitary normalizer of $\fD$.

\begin{prop}\label{P:invD}
Suppose that $\Fth$ is aperiodic and $\lambda_{(U,V)}\in\y$.  Let $W:=U\yp(V)(=V\ypp(U))$.
Then the following hold true.
\begin{itemize}
\item[(i)] The fixed point algebra for $\{\lambda_{(U,V)}:U,V\in\U(\fD)\}$ is $\fD$.
\item[(ii)] $\{\lambda_{(U,V)}: U,V\in\U(\fD)\}=\{\lambda\in\aut: \lambda|_\fD=\id\}$.
\item[(iii)] $\lambda_{(U,V)}(\fD)\subseteq\fD\Longleftrightarrow W\in N(\fD)$.
\\
If, in addition, $\lambda_{(U,V)}\in\aut$, then
$$\lambda_{(U,V)}(\fD)=\fD\Longleftrightarrow W\in N(\fD).$$
\item[(iv)] If $U,V\in\fD$, then
$\lambda_{(U,V)}\in\aut\Longleftrightarrow \lambda_{(U,V)}(\fD)=\fD$.
\end{itemize}
\end{prop}

\begin{proof}
(i)  To simply our writing, we use $\F ix$ to denote the fixed point algebra of
$\{\lambda_{(U,V)}:U,V\in\fD\}$.

In order to check  $\fD\subseteq\F ix$, it is sufficient to check that all generators of $\fD$ are in
$\F ix$. Represent a generator $X$ of $\fD$ as
$$X=s_{u_1\cdots u_k}s_{u_1\cdots u_k}^*,$$
where $u_i\in\Fth$ with $d(u_i)=(1,1)$ ($i=1,...,k$). 
Then, since both
$W$ and $s_{u_\ell\cdots u_k}s_{u_\ell\cdots u_k}^*$ ($1\le \ell \le k$) are in $\fD$,
they commute. Hence
we have
\begin{alignat*}{2}
\lambda_{(U,V)}(X)
&=\lambda(s_{u_1})\cdots \lambda(s_{u_k})\lambda(s_{u_k})^*\cdots \lambda(s_{u_1})^*\\
&=Ws_{u_1}\cdots Ws_{u_k}s_{u_k}^* W^*\cdots s_{u_1}^*W^*\\
&=Ws_{u_1}\cdots s_{u_k}s_{u_k}^*WW^*\cdots s_{u_1}^*W^*\\
&=\cdots \\
&=X.
\end{alignat*}
Thus $\fD\subseteq\F ix$.

We now prove $\F ix\subseteq\fD$. For any $D\in\U(\fD)$, by Lemma \ref{L:mult} we have
$\Ad(D)\in\{\lambda_{(U,V)}:U,V\in\U(\fD)\}$. So if $X\in\F ix$, then $\Ad(D)(X)=X$.
That is, $DXD^*=X$ for all $D\in\U(\fD)$. So we have $X\in\fD$ as $\fD$ is a masa of $\O_\theta$ by Lemma \ref{L:masa}.
This takes care of (i).

(ii) Assume that $\lambda_{(U,V)}|_\fD=\id$. For any $D\in\fD$,
we have $s_{e_i}Ds_{e_i}^*=\lambda_{(U,V)}(s_{e_i}Ds_{e_i}^*)=Us_{e_i} D s_{e_i}^*U^*$
($i=1,...,m$).
So $U\in\fD$. Similarly, we have $V\in \fD$. Hence
$$\{\lambda_{(U,V)}: U,V\in\U(\fD)\}\supseteq \{\lambda\in\y: \lambda|_\fD=\id\}.$$
The inclusion $\subseteq$ is directly from (i). This proves (ii).
Clearly, by Lemma \ref{L:mult} (ii) every element in the above set is an automorphism.

(iii) Set $\lambda=\lambda_{(U,V)}$.
Assume first $W\in N(\fD)$. From \eqref{E:limInv} we have
$\lambda(\fD)\subseteq\fD$.
Now suppose $\lambda(\fD)\subseteq\fD$. For $1\le i\le m, 1\le j\le n$, we have
$\fD\supseteq\lambda(s_{e_if_j}\fD s_{e_if_j}^*)=Ws_{e_if_j}\fD s_{e_if_j}^*W^*$. Thus $W\in N(\fD)$.

To finish the proof of (iii), it suffices to show that if $\lambda\in\aut$, then 
 $\lambda(\fD)\subseteq\fD$ actually implies $\lambda(\fD)=\fD$. 
 In the sequel, we show $\lambda^{-1}(\fD)\subseteq\fD$.
Arbitrarily choose $D_1,D_2\in\fD$. As $\lambda(\fD)\subseteq\fD$, we have 
\[
\lambda(\lambda^{-1}(D_1)D_2)=D_1\lambda(D_2)=\lambda(D_2)D_1=\lambda(D_2\lambda^{-1}(D_1)).
\]
As $\lambda$ is an automorphism, we have
$\lambda^{-1}(D_1)D_2=D_2\lambda^{-1}(D_1)$.
This implies
$\lambda^{-1}(D_1)\in\fD'$. As $\fD$ is a masa in $\O_\theta$ by Lemma \ref{L:masa},
it follows that $\lambda^{-1}(D_1)\in\fD$. Therefore,
$\lambda^{-1}(\fD)\subseteq\fD$ from the arbitrariness of $D_1$.

(iv) Clearly $W\in N(\fD)$ as $U,V\in\fD$. From (iii), $\lambda_{(U,V)}(\fD)\subseteq\fD$.
Suppose $\lambda_{(U,V)}\in\aut$. Then, from the proof of (iii), we see that 
 $\lambda_{(U,V)}(\fD)=\fD$. The proof of the other direction can be easily adapted from 
Proposition \ref{P:invF} (ii), and so omitted here. 
\end{proof}

From the proofs of Proposition \ref{P:invD}, we remark the following. First of all, 
if $(U,V)\in\U(\fD)\times \U(\fD)$ determines 
an endomorphism, then it is automatically an automorphism. 
Secondly, in (iii), we indeed have 
$W\in N(\fD)\Leftrightarrow\lambda_{(U,V)}(\fD)=\fD\Leftrightarrow \lambda_{(U,V)}(\fD)\subseteq\fD.$

As a consequence of Proposition \ref{P:invD} (i) and Lemma \ref{L:mult}, we get

\begin{cor}
Suppose $\Fth$ is aperiodic. Then 
 $\{\lambda_{(U,V)}\in\End(\O_\theta):U,V\in\U(\fD)\}$ is a maximal abelian subgroup of $\aut$.
\end{cor}

\begin{proof}
By Proposition \ref{P:invD} (ii), every element in 
$\{\lambda_{(U,V)}\in\End(\O_\theta):U,V\in\U(\fD)\}$ 
is an automorphism.
It now suffices to notice that, for $D\in\U(\fD)$, the identity
$\lambda_{(U,V)}\Ad(D)=\Ad(D)\lambda_{(U,V)}$ implies that 
$D^*\lambda_{(U,V)}(D)\in\Z(\O_\theta)=\bC I$.
\end{proof}

\subsection{Unitarily implemented automorphisms}
\label{SubS:uni}

In \cite{Voi}, Voiculescu constructed a family of unitarily implemented automorphisms
of the Cuntz algebras $\O_n$ from a subgroup $\U(n,1)$ of the general linear group $GL_n(\bC)$.
This result plays a very important role in many places. See, e.g., \cite{DP0, P1}.
Our original main purpose was to ``naturally'' generalize this result to 2-graph algebras.
But, as we have mentioned in Section \ref{S:End},
for a given pair of unitaries $(U,V)$ of $\O_\theta$, in practice, it is hard to check 
if $(U,V)$ determines an endomorphism
because of the twisted property.
So to know if it gives an automorphism becomes a much more challenging  task.
Thus, in this direction, so far we are only able to generalize the above result in the case of $\theta=\id$. 

In order to state our results, we first need some notation.
Following \cite{CE}, let
$ J=\begin{bmatrix} -1&0\\ 0& I_n\end{bmatrix}$
and
$$\U(n,1)=\left\{A=\begin{bmatrix}a_0 & h_1^*\\ h_2 & A_1\end{bmatrix}\in GL_{n+1}(\bC): A^*JA=J\right\}.$$
Here $a_0\in\bC$, $A_1$ is an $n\times n$ matrix, and $h_1,h_2$ are column vectors in $\bC^n$.
It is well-known that for each $A\in \U(n,1)$, there is a unitary $U_A\in\O_n$ 
determined by $A$, whose formula can be found in \cite{CE, Voi}.

We are now ready to give a family of unitarily  implemented automorphisms of $\O_{\id}$ 
constructed from $\U(m,1)\times \U(n,1)$.

\begin{prop}
Every pair $(A,B)\in \U(m,1)\times \U(n,1)$ determines a unitarily implemented
automorphism of $\O_{\id}$. Furthermore,
there is an action $\alpha$ of
$\U(m,1)\times \U(n,1)$ on $\O_{\id}$ given by
$$\alpha(A,B)(s_{e_i})=U_A s_{e_i},\ \alpha(A,B)(s_{f_j})=V_Bs_{f_j}$$
for all $1\le i\le m, 1\le j\le n$.
Here, $U_A$ (resp. $V_B$) are unitaries determined by $A$ (resp. $B$)
in the Cuntz algebras 
$\O_m=\ca(s_{e_1},...,s_{e_m})$ (resp. $\O_n=\ca(s_{f_1},...,s_{f_n})$).
\end{prop}

\begin{proof}
\footnote{
This proof is due to the referee. It is easier and shorter than the original one. 
}
Let $(A,B)\in \U(m,1)\times \U(n,1)$. 
As in \cite{Cun2}, let $\lambda_{U_A}\in{\rm End}(\O_m)$ and $\lambda_{V_B}\in{\rm End}(\O_n)$
denote the endomorphisms determined by $U_A$ and $V_B$, respectively.
Then, from \cite[2.9]{Voi}, $\lambda_{U_A}$ and $\lambda_{V_B}$
are actually unitarily implemented automorphisms of $\O_m$ and $\O_n$,
respectively.
Since $\theta=\id$,  from Example \ref{Eg:id},  we have 
$\lambda_{(U_A,V_B)}\in{\rm Aut}(\O_\id)$ and 
$\lambda_{(U_A,V_B)}=\lambda_{U_A}\otimes\lambda_{V_B}$,
Therefore $\lambda_{(U_A,V_B)}$ is also unitarily implemented.

The mapping $\alpha$ given in the proposition is an action because 
simple calculations yield 
$$
U_{A_2A_1}=\alpha_{(A_2,B_2)}(U_{A_1})U_{A_2},\quad 
V_{B_2B_1}=\alpha_{(A_2,B_2)}(V_{B_1})V_{B_2}
$$
for all $A_1,A_2\in\U(m,1)$ and $B_1,B_2\in\U(n,1)$. 
\end{proof}

\section{Modular theory of 2-graph algebras}\label{S:modular}

Recall that $\Phi=\int_{\bT^2}\gamma_{\bt}d\bt$ is the faithful conditional expectation of $\O_\theta$
onto the $(mn)^\infty$-UHF algebra $\fF$. Let $\tau$ be the unique faithful normalized trace on $\fF$.
Define $\omega=\tau\Phi.$
Then $\omega$ is a faithful state on $\O_\theta$. Also notice that $\omega\gamma_\bt=\omega$
($\bt\in \bT^2$), i.e.,
$\omega$ is invariant under the gauge automorphisms $\gamma_\bt$ of $\O_\theta$.

Let $L^2(\O_\theta)$ be the GNS Hilbert space determined by the state $\omega$. So the inner
product on $\O_\theta$ is defined by $\langle A| B\rangle=\omega(A^*B)$ for all $A,B\in\O_\theta$.
(Notice that the inner product here is linear in the second variable.)
Let $A\in\O_\theta$ and denote
the left action of $A$ by $\pi(A)$, that is, $\pi(A)B=AB$ for all $B\in \O_\theta$.
Let ${\O_{\theta c}}$ denote the algebra as the finite linear span of the generators of $\O_\theta$:
${\O_\theta}_c=\spn\{s_us_v^*:u,v\in\Fth\}$.

For brevity, in what follows, let $\bn:=(m,n)$, where $m,n$ are the numbers of generators of $\Fth$ of degree 
$(1,0)$ and $(0,1)$, respectively. 

The first lemma below gives some identities on the tracial state $\tau$ on $\fF$ and generalizes 
\cite[Lemma 3.1]{CPR}.

\begin{lem}\label{L:tronF}
Suppose $u,v\in\Fth$ with $d(u)=d(v)$. Then
$$
\tau(s_u X s_v^*)
=\delta_{u,v}\bn^{-d(u)}\tau(X)
\quad \text{for all}\quad X\in\fF.
$$
Here, as usual, $\delta_{u,v}=1$ if $u=v$; 0, otherwise.
\end{lem}

\begin{proof} Let $X\in\fF$.
Since $d(u)=d(v)$, we have
\begin{align}\label{E:omega}
\tau(s_u X s_v^*)
=\tau(s_u X s_u^*s_u s_v^*)
=\tau(s_u s_v^* s_u X s_u^*)
=\delta_{u,v}\tau(s_u X s_u^*).
\end{align}
In particular, we have
$$\tau(s_u s_v^*)=\delta_{u,v}\tau(s_us_u^*).$$
Making use of the defect free property, we have
$$
\sum_{d(u)=d(v)}\tau\left(s_u  s_u^*\right)
 =\tau\left(\sum_{d(u)=d(v)}s_u  s_u^*\right)
=\tau(I)
=1.
$$
Since there are only $\bn^{d(v)}$ such $u$'s with $d(u)=d(v)$,
we obtain
$$\tau(s_u s_u^*)=\bn^{-d(u)}.$$
Therefore
\begin{eqnarray}
\label{E:trace}
\tau(s_u s_v^*)=\delta_{u,v}\bn^{-d(u)}.
\end{eqnarray}

Now we claim that for any $(p,q)\in\bZ_+^2$, we have
\begin{alignat*}{2}
\tau(X)=\sum_{d(v)=(p,q)}\tau(s_v X  s_v^*)\quad \text{for all}\quad X\in \fF.
\end{alignat*}
To this end, it suffices to check it for the generators
$X=s_{w_1}s_{w_2}^*$ of $\fF$.
But then from \eqref{E:trace} it follows that
\begin{alignat*}{2}
\sum_{d(v)=(p,q)}\tau(s_v X s_v^*)
&=\sum_{v}\tau(s_v s_{w_1}s_{w_2}^*s_v^*)\\
&=\sum_{v}\delta_{w_1,w_2}\bn^{-d(w_1)-d(v)}\quad (\mbox{by}\ \eqref{E:trace})\\
&=\delta_{w_1,w_2}\bn^{-d(w_1)}\\
&=\tau(X) \quad (\mbox{by}\ \eqref{E:trace}).
\end{alignat*}

On the other hand, for any $X\in\fF$,
similar to the proof of \eqref{E:omega} we have
\begin{alignat*}{2}
\tau(s_u X  s_u^*)
&=\tau(s_u X s_v^*s_v  s_u^*)\\
&=\tau(s_v  s_u^* s_u X s_v^*)\\
&=\tau(s_v X s_v^*).
\end{alignat*}
Thus
$$
\tau(s_u X s_u^*)
=\bn^{-d(u)}\tau(X)\quad  \mbox{for all}\quad X\in \fF.
$$
Combining this identity with \eqref{E:omega} completes the proof of the lemma.
\end{proof}

\bigskip

We now begin to give the modular objects in the celebrated Tomita-Takesaki modular theory.  
Define an operator $S$ on ${\O_\theta}_c\subset L^2(\O_\theta)$ by
$$S(A)=A^*\quad \mbox{for all}\quad A\in {\O_\theta}_c.$$
Clearly, $S$ is anti-linear.
Define another anti-linear operator $F$ on ${\O_\theta}_c$, which acts on
generators by
$$
F(s_us_v^*)
=\bn^{d(u)-d(v)}s_vs_u^*,
$$
and then extend it anti-linearly.

We shall show that $F$ is indeed the adjoint of $S$.
The key step in its proof is to make full use of the close relations between $\omega$ and $S$, $F$.
The basic idea behind here
is to convert the computations involved with
$\langle S(A)|B\rangle$ and $\langle F(B)|A\rangle$ to those related to
$\omega$. But then
$\omega|_{\fF}=\tau$ is a trace, and so we can invoke the commutativity of $\tau$ and apply
Lemma \ref{L:tronF}. If we use this approach to Cuntz algebras, it seems that the proof
here is more unified than that in \cite{CPR}.
More importantly,
if one applies the approach in \cite{CPR} directly, it seems that he/she could only deal with a very
special class of 2-graph algebras, i.e., those with the identity relation
($\theta=\id$).

\begin{lem}\label{L:FadjS}
Let $S,F$ be defined as above.
Then $F$ is the adjoint of $S$: $\langle S(A)|B\rangle=\langle F(B)|A\rangle$
for all $A,B\in{\O_\theta}_c$.
\end{lem}

\begin{proof}
It suffices to show that $\langle S(A)|B\rangle=\langle F(B)|A\rangle$
holds for all generators $A,B$ of ${\O_\theta}_c$.
So we let
$A=s_{u_1}s_{v_1}^*$ and 
$B=s_{u_2}s_{v_2}^*$.
From the definition of the degree map $d$ for the generators of $\O_\theta$,
we have $d(A)=d(u_1)-d(v_1)$ and $d(B)=d(u_2)-d(v_2)$. 

First observe from the definitions of $S, F$ that
\begin{eqnarray}
\label{E:Sw}
\langle S(A)|B\rangle=\langle A^*|B\rangle=\omega(AB)
\end{eqnarray}
and
\begin{eqnarray}
\label{E:Fw}
\langle F(B)|A\rangle
=\bn^{d(B)}\langle B^*|A\rangle
=\bn^{d(B)}\omega(BA).
\end{eqnarray}

If $d(A)+d(B)\ne 0$, then $AB,BA$ are either $0$, or not in $\fF$.
So \eqref{E:Sw} and \eqref{E:Fw} implies
$\langle S(A)|B\rangle=0=\langle F(B)|A\rangle$.
We are done.
Thus we now suppose that $d(A)+d(B)=0$ in what follows.

\smallskip
\textit{Case 1.}
$d(A)=(-s,-t)$ with $s,t\ge 0$. Then $d(B)=(s,t)$. So, making full use of commutation
relations of $\Fth$, we can rewrite $A,B$ as
\begin{align*}
A=A' s_u^*,\quad 
B =s_vB',
\end{align*}
for some $u,v\in\Fth$ with $d(u)=d(v)=(s,t)$ and some generators $A',B'$ in $\fF$.

Clearly, we now have
\begin{alignat}{4}
\label{E:wAB}
\omega(AB)
&=\omega(A' s_{u}^*s_vB')
=\delta_{u,v}\omega(A'B')
=\delta_{u,v}\tau(A'B').
\end{alignat}
On the other hand,
it follows from Lemma \ref{L:tronF} that
\begin{alignat}{4}\nonumber
\omega(BA)
&=\omega(s_vB'A's_u^*)\\\nonumber
&=\delta_{u,v}\bn^{-d(u)}\omega(B'A')\\\label{E:wBA}
&=\delta_{u,v}m^{-s}n^{-t}\tau(B'A').
\end{alignat}
As $\tau$ is a trace on $\fF$,  from \eqref{E:Sw}, \eqref{E:Fw}, \eqref{E:wAB} and \eqref{E:wBA} we proved
$\langle S(A)|B\rangle=\langle F(B)|A\rangle$.

\smallskip
\textit{Case 2.} $d(A)=(s,-t)$ with $s,t\ge 0$. Then $d(B)=(-s,t)$. As above, we
rewrite $A,B$ as
\begin{align*}
A=s_{u_1}A's_{v_1}^*, \quad 
B=s_{v_2}B's_{u_2}^*,
\end{align*}
for some $u_1,u_2,v_1,v_2\in\Fth$ and  some generators $A',B'\in\fF$ with 
$$d(u_1)=d(u_2)=(s,0),\quad d(v_1)=d(v_2)=(0,t).$$
We have from Lemma \ref{L:tronF} that
\begin{alignat*}{2}
\omega(AB)
&=\omega(s_{u_1}A's_{v_1}^*s_{v_2}B's_{u_2}^*)\\
&=\delta_{v_1,v_2}\omega(s_{u_1}A'B's_{u_2}^*)\\
&=\delta_{v_1,v_2} \delta_{u_1,u_2}\bn^{-d(u_1)} \omega(A'B')\\
&=\delta_{u_1,u_2} \delta_{v_1,v_2}m^{-s} \tau(A'B'),
\tag{$9'$}
\label{(9)'}
\end{alignat*}
and
\begin{alignat*}{2}
\omega(BA)
&=\omega(s_{v_2}B's_{u_2}^* s_{u_1}A's_{v_1}^* )\\
&=\delta_{u_1,u_2}\omega(s_{v_2}B'A's_{v_1}^* )\\
&=\delta_{u_1,u_2}\delta_{v_1,v_2}\bn^{-d(v_1)} \omega(B'A')\\
&=\delta_{u_1,u_2} \delta_{v_1,v_2}n^{-t} \tau(B'A').
\tag{$10'$}
\label{(10)'}
\end{alignat*}
It follows from \eqref{E:Sw}, \eqref{E:Fw}, \eqref{(9)'} and \eqref{(10)'} that
$\langle S(A)|B\rangle=\langle F(B)|A\rangle$.

\smallskip
\textit{Case 3.} $d(A)=(-s,-t)$ with $s,t\ge 0$. The proof is completely similar to Case 1.

\textit{Case 4.} $d(A)=(s,-t)$ with $s,t\ge 0$. The proof is completely similar to Case 2.

Therefore, $F$ is the adjoint of $S$.
\end{proof}

From Lemma \ref{L:FadjS}, we particularly obtained that both $F$ and $S$
are closable (\cite[Theorem VIII.1]{ReeSim}). By abusing notation, 
we still use $F,S$ to denote their corresponding closures.
The following is well-known:
$S$ and $F$ have polar decompositions
\begin{eqnarray*}
S=J\triangle^{\frac{1}{2}}=\triangle^{-\frac{1}{2}}J,\quad
F=J\triangle^{-\frac{1}{2}}=\triangle^{\frac{1}{2}}J,
\end{eqnarray*}
where  $\triangle=FS$, and $J$ is an anti-unitary operator with $J=J^*$ and $J^2=I$.
So far we have obtained all modular objects.

Moreover, by Lemma \ref{L:FadjS} we have
$$
J(s_us_v^*)
=\bn^{\frac{d(u)-d(v)}{2}}s_vs_u^*,
$$
and
$$
\triangle^z(s_us_v^*)
=\bn^{z(d(v)-d(u))}s_us_v^* \quad(z\in\bC).
$$
Here $m^z=\exp(z\ln m)$.

We are now in a position to prove the first main result in this section.

\begin{thm}
\label{T:modalg}
The algebra ${\O_\theta}_c$ with the inner product $\langle\cdot | \cdot \rangle$:
$\langle A|B\rangle=\omega(A^*B)$,
is a modular Hilbert algebra.
\end{thm}

\begin{proof}
The proof can now be easily adapted from \cite[Lemma 3.2]{CPR}.
In order to check all axioms of a modular Hilbert
algebra (which is called a Tomita algebra in \cite{CPR}),
the only thing that is not very obvious here is  the fact that
every $\triangle^z$ is multiplicative on ${\O_\theta}_c$. We will prove this below.

Arbitrarily choose two generators of ${\O_\theta}_c$:
$
X=s_{u_1}s_{v_1}^*
$
and 
$
Y=s_{u_2}s_{v_2}^*.
$
Let $d(u_2)\bigvee d(v_1)=(p,q)$.
Then
$$
s_{v_1}^*s_{u_2}
=\sum s_{w_1}s_{w_2}^*,
$$
where the sum is over all $w_1,w_2\in\Fth$
such that 
\begin{align}\label{pq}
v_1w_1=u_2w_2 \quad\text{and}\quad d(w_1)+d(v_1)=(p,q)
\end{align}
 (cf. \cite{KumPask, Raeburn}).
We now have
\begin{alignat*}{3}
\triangle^z(XY)
&=\triangle^z\left(s_{u_1}\; \sum_{w_1,w_2}s_{w_1}s_{w_2}^* \; s_{v_2}^*\right)\\
&=\sum_{w_1,w_2}  \bn^{z(k_1,\;\ell_1)} 
      s_{u_1w_1}s_{v_2w_2}^*\\
&=\sum_{w_1,w_2}  \bn^{z(k_2,\;\ell_2)} 
      s_{u_1w_1}s_{v_2w_2}^*\quad (\text{by}\ \eqref{pq})\\
&= \bn^{z(k_2,\;\ell_2)} \sum_{w_1,w_2} 
      s_{u_1w_1}s_{v_2w_2}^*\\
&=\bn^{z(k_2,\;\ell_2)} XY\\
&=\triangle^z(X)\triangle^z(Y),
\end{alignat*}
where
\begin{align*}
(k_1,\ell_1)&:= d(v_2)+d(w_2)-d(u_1)-d(w_1),\\
(k_2,\ell_2)&:= d(v_1)+d(v_2)-d(u_1)-d(u_2).
\end{align*}
Therefore $\triangle^z$ is multiplicative on ${\O_\theta}_c$.
\end{proof}

The above operators $S$ (resp. $J$, $\triangle$) are called the \textit{Tomita operator}
(resp. \textit{modular conjugation}, \textit{modular operator}) of $\O_{\theta c}$.
Let $\pi(\O_\theta)''$ be the von Neumann algebra generated by the GNS representation
of the state $\omega$. Then $\pi(\O_\theta)''$ is nothing but the left von Neumann algebra of 
$\O_{\theta c}$ (\cite{Take}).
The celebrated Tomita-Takesaki modular
theory says that
$$
\triangle^{it}\pi(\O_\theta)''\triangle^{-it}=\pi(\O_\theta)''\ (t\in\bR),
\quad
J\pi(\O_\theta)''J=\pi(\O_\theta)'.
$$
Refer to \cite{Str, Take} for more information on the Tomita-Takesaki modular theory.

Let $\sigma_{z}(\pi(X))=\triangle^{iz}\pi(X)\triangle^{-iz}$ for all 
$z\in\bC$ and $X\in\O_{\theta c}$. 
We now give the formula of the modular automorphisms $\sigma_{t}$ $(t\in\bR)$ of $\pi(\O_\theta)''$
on generators, and 
show that $\omega$ is a $\sigma$-KMS state. Refer to \cite[Chapter 5]{BraRob}
and \cite[Section 13]{Take} for KMS-states.

\begin{prop}\label{P:tri}
Let $\omega$ be the state given at the beginning of this section, $\triangle$
the modular operator of $\O_{\theta c}$, and $\sigma_z$ $(z\in\bC)$ the operator 
defined as above. Then
\begin{itemize}
\item[(i)]
the group of modular automorphisms $\sigma_t$ ($t\in\bR$)
of the von Neumann algebra
$\pi(\O_\theta)''$ on $L^2(\O_\theta)$ is given on the generators by
\begin{align*}
\sigma_t(\pi(s_us_v^*))
:=\triangle^{it}\pi(s_us_v^*)\triangle^{-it}
=\bn^{it(d(v)-d(u))}\pi(s_us_v^*);
\end{align*}
\item[(ii)]
$\omega$ is a $\sigma$-KMS state over $\pi(\O_\theta)''$:
$$
\omega(AB)=\omega(\sigma_i(B)A) \quad \text{for all}\quad A,B\in \pi({\O_\theta}_c);
$$
\item[(iii)]
$\omega$ is the unique $\sigma$-KMS state over $\pi(\O_\theta)''$, provided  $\frac{\ln m}{\ln n}\not\in\bQ$.
\end{itemize}
\end{prop}

\begin{proof}
The proofs of (i) and (ii) below are borrowed from
\cite[Lemma 3.3 and its remarks]{CPR}.

(i) From the proof of Theorem \ref{T:modalg}, we know that
$\triangle^z:{\O_\theta}_c\to{\O_\theta}_c$ is an
algebra homomorphism. This implies that
$$\pi(\triangle^z(A))=\triangle^z\pi(A)\triangle^{-z} \quad \mbox{for all}\quad A\in{\O_\theta}_c,\ z\in\bC.$$
The rest of the proof of (i) is done by direct computation.

(ii)  It is proved by the following calculations:
\begin{alignat*}{2}
\omega(AB)&=\langle A^*|B\rangle=\langle S(A)|B\rangle=\langle F(B)|A\rangle=\omega(SF(B)A)\\
                     &=\omega(\triangle^{-1}(B)A)
                      =\omega(\sigma_i(B)A).
\end{alignat*}

(iii) In the sequel, we naturally identify the C*-algebra $\pi(\O_\theta)$ with the C*-algebra $\O_\theta$
(not as a subspace of $L^2(\O_\theta)$)
as $\pi$ is faithful.
A straightforward calculation yields the following relation:
$$
\sigma_t=\gamma_{(m^{-it},\; n^{-it})}\quad \text{for all}\quad t\in\bR.$$
As
$\frac{\ln m}{\ln n}\not\in\bQ$,
by Kronecker's Theorem,
the set $\{(m^{-it}, n^{-it}): t\in\bR\}$ is dense in $\bT^2$.
Thus  the modular automorphsims $\{\sigma_t:t\in\bR\}$
determine the gauge automorphisms $\{\gamma_\bt:\bt\in \bT^2\}$.

Now we can use a similar argument of \cite[Theorem 2]{OlePed}
(also cf. \cite[Chapter 5]{BraRob}) to prove
the uniqueness of $\omega$. We only sketch it here. Suppose $\omega'$ is also a KMS state for $\sigma$
at value $\beta$.

We first show that $\beta$ has to be finite. To the contrary, suppose that 
$\beta=\infty$ (or $-\infty$). However, from (i)
one can see that the functions
$$
t\mapsto \omega'(\pi(s_{e_i})^*\sigma_t(\pi(s_{e_i})))
=m^{-it}\omega'(\pi(s_{e_i})^*\pi(s_{e_i}))=m^{-it}\omega'(I)=m^{-it},
$$
$$({\rm or}\quad
t\mapsto \omega'(\pi(s_{e_i})\sigma_t(\pi(s_{e_i})^*))=m^{it}\omega'(\pi(s_{e_i}s_{e_i}^*)))
$$
do not have bounded analytic extensions to the upper (or lower) half-planes.
Here we used the simple fact that 
$\omega'(\pi(s_{e_i}s_{e_i}^*))\ne 0$ because of the defect free property.
Now from \cite[Proposition 5.3.19]{BraRob} and its immediately preceding remark,
we get a contradiction. 
Therefore, $\beta$ is finite.

Since $\omega'$ is a $(\sigma,\beta)$-KMS state, $\omega'$ is invariant under $\sigma_t$, namely,
$\omega'\sigma_t= \omega'$ ($t\in\bR$).
Hence from the relations between $\sigma$ and $\gamma$ given above, we obtain
$\omega'\gamma_\bt= \omega'$ for all $\bt\in\bT^2$. Hence, $\omega' \Phi=\omega'$.
On the other hand, from the KMS condition (i.e.,
$\omega'(AB)=\omega'(\sigma_{i\beta}(B)A)$),
we have that $\omega'|_{\fF}$ 
is a normalized trace on $\fF$. The uniqueness of the normalized trace on $\fF$ concludes
$\omega'\equiv \omega$. Now from (ii), we also have $\beta=1$. 
Therefore, $\omega$ is the unique $\sigma$-KMS state over $\pi(\O_\theta)''$.
\end{proof}

\section{Some remarks on the classification of $\pi(\O_\theta)''$}\label{S:class}

In this short section, we begin with the following observation.

\begin{lem}\label{L:Ffactor}
$\pi(\fF)''$ is a II$_1$ factor.
\end{lem}

\begin{proof}
It is known that $\pi$ and the GNS representation of $\omega|_{\fF}$
(the restriction $\omega|_{\fF}$ of $\omega$ to
$\fF$, and so $\omega|_{\fF}=\tau$) are quasi-equivalent.
The lemma now follows from \cite[Theorem 2.5]{Pows}.
\end{proof}

Recall that if $\fM$ is a von Neumann algebra, the \textit{Connes invariant $\S(\fM)$} is the intersection over all faithful
normal states of the spectra of their corresponding modular operators \cite{Jones}.
Connes classified type III factors as follows.
A factor $\fM$ is said to be of
\begin{alignat*}{3}
type\ III_0 \quad \text{if} & \quad \S(\fM)=\{0, 1\};\\
type\ III_\lambda \quad \text{if} &\quad  \S(\fM)=\{0, \lambda^n:n\in \bZ\} \quad (0<\lambda<1);\\
type\ III_1 \quad \text{if} &\quad  \S(\fM)=\{0\}\cup\bR^+.
\end{alignat*}

We are now able to obtain the following result on partially classifying the von Neumann algebra $\pi(\O_\theta)''$.

\begin{cor}\label{C:class}
 If $\frac{\ln m}{\ln n}\not\in\bQ$, then
$\pi(\O_\theta)''$ is an AFD factor of type III$_1$.
\end{cor}

\begin{proof}
Since $\O_\theta$ is amenable (\cite{KumPask}), it is known that $\pi(\O_\theta)''$ is AFD.
Furthermore, by Proposition \ref{P:tri} (iii) and \cite[Theorem 5.3.30]{BraRob},  we obtain that
$\pi(\O_\theta)''$ is a factor.

Also, from Proposition \ref{P:tri}, one can show that the fixed point algebra of $\sigma$ is $\pi(\fF)''$
as $\frac{\ln m}{\ln n}\not\in\bQ$.
It follows from Lemma \ref{L:Ffactor} and \cite[Section 28.3]{Str}  that the Connes spectrum coincides
with the spectrum of the modular operator. That is, $\S(\pi(\O_\theta)'')=\rm{Sp}(\triangle)$.
But from the definition of $\triangle$, it
is easy to see that
$$\rm{Sp}(\triangle)=\overline{\{m^an^b:a,b\in\bZ\}}=[0,\infty)$$
as $\frac{\ln m}{\ln n}\not\in\bQ$.
So $\pi(\O_\theta)''$ is of type III$_1$.
\end{proof}

\begin{rem}
Notice that if $\frac{\ln m}{\ln n}\in\bQ$, i.e., $m^a=n^b$ for some $a,b\in\bN$ with $\rm{gcd}(a,b)=1$,
then $\pi(\fF)''$ is a \textit{proper} subalgebra of $\pi(\O_\theta)''_\sigma$,
the fixed point algebra of the modular
automorphisms $\sigma_t$ ($t\in\bR$).  Indeed, one can check that
\begin{align*}
\pi(\O_\theta)''_\sigma
&=\{\pi(\fF), \pi(s_{e_{u_1}}s_{f_{v_1}}^*), \pi(s_{f_{v_2}}s_{e_{u_2}}^*):
                    (|u_{i}|,|v_{i}|)=k_{i}(a, b), k_{i}\in\bN\}'' \\
&\supsetneq \pi(\fF)''.
\end{align*}
\end{rem}

\begin{rem}
 Since $\frac{\ln m}{\ln n}\not\in\bQ$, the 2-graph $\Fth$ in Corollary \ref{C:class} is aperiodic.
Clearly, if $\Fth$ is periodic, then $\pi(\O_\theta)''$ is not a factor
(cf. \cite{DY}). A natural question is if the converse is true:  If $\Fth$ is aperiodic, is $\pi(\O_\theta)''$ a factor?
Moreover,
by Corollary \ref{C:class}, we obtain the type of $\pi(\O_\theta)''$ when $\frac{\ln m}{\ln n}\not\in\bQ$.
If the converse is true, a further question arises: What is the type 
of $\pi(\O_\theta)''$ if $\frac{\ln m}{\ln n}\in\bQ$ and $\Fth$ is aperiodic?
Also, it would be interesting to study the index theory of
endomorphisms of $\pi(\O_\theta)''$.
We will leave those as future research.
\end{rem}

\medskip
\noindent
\textbf{Acknowledgements.}
I would like to thank Professor Kenneth R. Davidson for his continuous encouragement and
helpful discussion. I also thank Professor W. Szyma\'nski for generously sending me his paper \cite{Szy}
and pointing out one error in the original proof of Proposition \ref{P:invF} (i).
I am very grateful to the anonymous referee for his/her very careful reading and 
providing many valuable comments which
greatly improve the presentation of the paper.

\end{document}